\documentclass[11pt]{amsart}
\usepackage{Cycle}
\usepackage[all]{xy}

\newcommand{\CartS}[1][dr]{\ar@{}[#1]|-*+{\square}}

\begin{document}

\title{The \'Etale Homology and The Cycle Maps in Adic Coefficients}
\author{Ting Li}
\address{Department of Mathematics\\
    Sichuan University\\
    Chengdu 610064\\
    P. R. China}
\email{moduli@scu.edu.cn}

\maketitle

\begin{abstract}
In this article, we define the $\ell$-adic homology for a morphism of schemes satisfying certain finiteness conditions. This homology has these functors similar to the Chow groups: proper push-forward, flat pull-back, base change, cap-product, etc. In particular on singular varieties, this kind of $\ell$-adic homology behaves much better that the classical $\ell$-adic cohomology. As an application, we give an much easier approach to construct the cycle maps for arbitrary algebraic schemes over fields of finite cohomology dimension. And we prove these cycle maps kill the algebraic equivalences and commute with the Chern action of locally free sheaves.
\end{abstract}
\plainfoot{Mathematics Subject Classification 2000: 14F20, 14C25}
\plainfoot{This work is supported by NSFC(10626036)}

\section*{Introduction}

The \'etale cohomology, especially the $\ell$-adic cohomology, is one of the most important tools of modern algebraic and arithmetic geometry, which allows us to construct a ¡°good¡± cohomology theory for varieties over fields of arbitrary characteristic. More specifically, people use the $\ell$-adic cohomology $\upH{\ast}(\etSite{X}, \mathbb{Z}_{\ell})$ to substitute for singular cohomology on varieties of arbitrary characteristic. On a nonsingular varieties, the cohomology $\upH{\ast}(\etSite{X}, \mathbb{Z}_{\ell})$ has very good properties and produces rich results. But on singular varieties or more generally on arbitrary schemes, the cohomology $\upH{\ast}(\etSite{X}, \mathbb{Z}_{\ell})$ behave not so good, and many important constructions and results are not valid. So on singular varieties, the \'etale homology is more suitable than the \'etale cohomology.

In this paper, we generalize the \'etale homology defined in \cite{MR0444667} in the following three facets. First, we define the \'etale homology in adic coefficients, which we call the $\ell$-adic homology. Second, our theory of $\ell$-adic homology is defined over schemes separated and of finite type over base schemes satisfying certain finiteness conditions, not just the algebraic schemes over separably closed fields as in \cite{MR0444667}. In particular, algebaic schemes over fields which are not necessarily separably closed, are considered by us. Since our theory is based on the adic formalism created by Ekedahl \cite{MR1106899}, the $\ell$-adic homology over base schemes of certain finiteness conditions shares almost the same good functorial properties, with that over separably closed base fields. Third, the $\ell$-adic homology groups $\rmH{\ast}(X, \mathscr{F})$ defined by us take value in arbitrary bounded complex $\mathscr{F}$, not just $\mathbb{Z}_{\ell}$, $\mathbb{Q}_{\ell}$ or $\mathbb{Z}/n\mathbb{Z}$ as in \cite{MR0444667}. And almost all functors and properties are preserved when extending to complexes.

In \S\ref{Sec1}, we briefly reiterate the category $\conD{\etSite{X}, R_{\smdot}}$ together with  the Grothendieck's six operations in \cite{MR1106899}. In \S\ref{Sec5}, we recite the properties of the functor $\rdF f^!$ and use the language of \cite{MR1106899} to rewrite the trace morphisms introduced in \cite[XVIII]{SGA4} and \cite[Cycle]{MR0463174}.

In \S\ref{Sec6}, we define the $\ell$-adic homology groups $\rmH{n}(X/Y, \mathscr{N})$ and $\mathbb{H}_n(X/Y, \mathscr{N})$ for a morphism $X \to Y$ of schemes satisfying certain finiteness conditions. These homology groups behave similarly in many facets to the bivariant Chow groups \mbox{$A^{-n}(X \to Y)$} defined in \cite[Ch.~17]{MR732620}. If $X$ is a $d$-dimensional nonsingular variety over a separably closed field $k$, then
\[ \mathbb{H}_{n}(X/k, \mathscr{N}) = \upH{2(d-n)}\bigl(\etSite{X}, \mathscr{N}(d-n)\bigr) \,. \]
We define two maps: the push-forward maps $f_{\ast}$ and the pull-back maps $f^{\ast}$, for the $\ell$-adic homology groups, which correspond to these maps on Chow groups $\cycA_{\ast}(X)$ defined in \cite[\S1.4 \& \S1.7]{MR732620}. We prove that the two maps $f_{\ast}$ and $f^{\ast}$ commute (see Theorem \ref{217P3}), which is essential to the construct various cycle maps basing on $\ell$-adic homology. Moreover we define the base change maps on the $\ell$-adic homology.

In \S\ref{Sec7}, we apply the $\ell$-adic homology in \S\ref{Sec6} to define the cycle map
\[ \clM_{X,\ell} \colon \cycA_{\ast}(X) \to  \mathbb{H}_{\ast}(X, \mathbb{Z}_{\ell}) \]
for arbitrary algebraic scheme $X$ over a field of finite cohomological dimension at $\ell$. We prove that the cycle map $\clM_{X,\ell}$ commutes with the push-forward map $f_{\ast}$ and the pull-back map $f^{\ast}$. And we prove the cycle maps kill the algebraic equivalence of algebraic cycles. In \S\ref{Sec8}, we prove that the cycle map $\clM_{X,\ell}$ commutes with the Chow action $c_i(\mathscr{E}) \cap \smdot[.6em]$ by locally free sheaves.

\medskip

{\bfseries Notation and Conventions.} A morphism $f \colon X \to Y$ of schemes is said to \emph{flat (resp.~smooth) of relative dimension $n$} if $f$ is flat (resp.~smooth) and all fibers of $f$ are $n$-equidimensional.

A morphism $f \colon X \to Y$ of Noetherian schemes is said to be \emph{compactifiable} if it factors as $f = \bar{f} \circ j$ where $j \colon X \hookrightarrow \overline{X}$ is an open immersion, and $\bar{f} \colon \overline{X} \to Y$ is a proper morphism. By \cite[Theorem 4.1]{MR2356346}, $f$ is compactifiable if and only if it is separated and of finite type.

An \emph{algebraic scheme} over a field $k$ is a scheme separated, of finite type over $k$. A \emph{variety} over $k$ is an integral algebraic scheme over $k$.

If $A$ is a Noetherian ring, we define $\fgD{A}$ to be the full subcategory of $\dD{A}$ consisting of complexes cohomologically finitely generated.

If $\mathscr{F}^{\smdot}$ is a complex of sheaves on $\etSite{X}$, we write $\mathscr{F}^{\smdot}\dangle{r} \coloneqq \mathscr{F}^{\smdot}(r)[2r]$ for each $r \in \mathbb{Z}$.

The notation $\coloneqq$ means being defined as; $\isoTo$ means isomorphism; and the notation $\square$ in commutative diagrams means Cartesian square.

\section{The \texorpdfstring{$\ell$}{l}-adic sheaves}\label{Sec1}

In this section, we briefly reiterate the theory of Ekedahl \cite{MR1106899} about the category $\conD{\etSite{X}, R_{\smdot}}$ together with  the Grothendieck's six operations. See also \cite{MR1963494} and \cite{YLMOp}.

Fix a prime number $\ell$, and let $R$ be the integral closure of $\mathbb{Z}_{\ell}$ in a finite extension field of $\mathbb{Q}_{\ell}$.

Let $X$ be a Noetherian scheme. We denote by $\mathbf{S}(\etSite{X}^{\mathbb{N}}, R_{\smdot})$ the abelian category of inverse systems
\[ \cdots \to \mathscr{F}_{n+1} \xrightarrow{p_n} \mathscr{F}_n \to \cdots \to \mathscr{F}_2 \xrightarrow{p_1} \mathscr{F}_1 \]
such that each $\mathscr{F}_n$ is a sheaf of $R_n$-modules on $\etSite{X}$. Set
\[ \dD{\etSite{X}^{\mathbb{N}}, R_{\smdot}} \coloneqq \dD(1){\mathbf{S}(\etSite{X}^{\mathbb{N}}, R_{\smdot})} \,; \]
and let $\conD{\etSite{X}^{\mathbb{N}}, R_{\smdot}}$ be the full subcategory of $\dD{\etSite{X}^{\mathbb{N}}, R_{\smdot}}$ consisting of complexes cohomologically AR-adic and constructible. Let $\conD{\etSite{X}, R_{\smdot}}$ be the quotient of $\conD{\etSite{X}^{\mathbb{N}}, R_{\smdot}}$ by inverting AR-quasi-isomorphisms.

If $f \colon X \to Y$ is a morphism of Noetherian schemes, then we have a triangulated functor
\[ f^{\ast} \colon \conD{\etSite{Y}, R_{\smdot}} \to \conD{\etSite{X}, R_{\smdot}} \,. \]

As to other five operations, we must add some restrictions on the underlying schemes. We consider the following condition (\dag) relate to a scheme $X$:
\begin{quotation}
(\dag) $X$ is Noetherian, quasi-excellent, of finite Krull dimension; $\ell$ is invertible on $X$ and $\mathrm{cd}_{\ell}(X) < \infty$.
\end{quotation}

From the Gabber's finitenes theorem for \'etale cohomology in \cite{OGab1}, we know the following facts
\begin{enumerate}
\item If $X$ satisfies (\dag), then any scheme of finite type over $X$ satisfies (\dag).
\item Let $R$ be a quasi-excellent, Henselian local ring with residue field $k$ such that $\mathrm{cd}_{\ell}(k) < \infty$. Then $\Spec R$ satisfies (\dag).
\item If $\ell \neq 2$, then the affine scheme $\Spec \mathbb{Z}[1/\ell]$ satisfies (\dag). (See \cite[X, 6.1]{SGA4})
\item If $f \colon X \to Y$ is a compactifiable morphism of schemes satisfying (\dag); then both $\rdF f_{\ast}$ and $\rdF f^!$ are of finite cohomological amplitude.
\end{enumerate}

In particular if $X$ is a scheme satisfying (\dag), then $\etSite{X}$ satisfies the condition A) in \cite{MR1106899}; thus we have two bi-triangulated functors
\begin{align*}
\smdot[.4em] \ttP[R] \smdot[.4em] & \colon \conDn(1){\etSite{X}, R_{\smdot}} \times \conDn(1){\etSite{X}, R_{\smdot}} \to \conDn(1){\etSite{X}, R_{\smdot}} \,, \\
\sRHom_R(\smdot[.6em], \, \smdot[.6em]) & \colon \conDn(1){\etSite{X}, R_{\smdot}}^{\mathrm{opp}} \times \conDp(1){\etSite{X}, R_{\smdot}} \to \conDp(1){\etSite{X}, R_{\smdot}} \,.
\end{align*}
And if $f \colon X \to Y$ is a compactifiable morphism of schemes satisfying (\dag), there are triangulated functors
\begin{align*}
\rdF f_{\ast} & \colon \conD{\etSite{X}, R_{\smdot}} \to \conD{\etSite{Y}, R_{\smdot}} \,, \\
\rdF f_! & \colon \conD{\etSite{X}, R_{\smdot}} \to \conD{\etSite{Y}, R_{\smdot}} \,, \\
\rdF f^! & \colon \conD{\etSite{Y}, R_{\smdot}} \to \conD{\etSite{X}, R_{\smdot}} \,.
\end{align*}

For each scheme $X$ satisfying (\dag), each object $\mathscr{F}$ in $\conD{\etSite{X}, R_{\smdot}}$, and each $n \in \mathbb{N}$, we define
\[ \upH{n}(\etSite{X}, \mathscr{F}) \coloneqq \Hom_{\conD{\etSite{X}, R_{\smdot}}}(R_X, \mathscr{F}[n]) \,. \]
Note that this definition is compatible with the continuous \'etale cohomology $\mathrm{H}_{\mathrm{cont}}^n(\etSite{X}, \mathscr{F})$ defined in \cite{MR929536}.

When we consider the schemes of finite type over a separably closed field, the following Theorem is essential.

\begin{Thm}\label{194P14}
The right derived functors of $(M_n) \mapsto \varprojlim M_n$ and the left derived functors of $M \mapsto (M \otimes_R R_n)$ define a natural equivalence of categories between $\conD{R_{\smdot}}$ and $\fgD{R}$.
\end{Thm}

\begin{proof}
See \cite[Proposition 2.2.8]{MR1963494}.
\end{proof}

Now we fix a separably closed field $k$. Note that $\conD(1){\etSite{(\Spec k)}, R_{\smdot}} = \conD{R_{\smdot}} = \fgD{R}$.

\begin{Nota}
Let $X$ be an algebraic scheme over $k$, $p \colon X \to \Spec k$ the structural morphism. Put
\begin{align*}
\rdF \varGamma(\etSite{X}, \smdot[.6em]) & \coloneqq \rdF p_{\ast} \colon \conD{\etSite{X}, R_{\smdot}} \to \fgD{R} \,, \\
\rdF \varGamma_!(\etSite{X}, \smdot[.6em]) & \coloneqq \rdF p_! \colon \conD{\etSite{X}, R_{\smdot}} \to \fgD{R} \,.
\end{align*}
Then for each $q \in \mathbb{Z}$, we have $\upH{q}(\etSite{X}, \smdot[.6em]) = \upH{q} \circ \rdF \varGamma(\etSite{X}, \smdot[.6em])$. And we define
\[ \upHc{q}(\etSite{X}, \smdot[.6em]) \coloneqq \upH{q} \circ \rdF \varGamma_!(\etSite{X}, \smdot[.6em]) \,. \]
\end{Nota}

\begin{Thm}[The K\"unneth Formula]\label{194T2}
Let $X$ and $Y$ be two algebraic schemes over $k$, $Z \coloneqq X \times_k Y$, $f \colon Z \to X$ and $g \colon Z \to Y$ the projections. Then for each $\mathscr{F} \in \conDn{\etSite{X}, R_{\smdot}}$ and $\mathscr{G} \in \conDn{\etSite{Y}, R_{\smdot}}$, there are two natural isomorphisms in $\fgDn{R}$\textup{:}
\begin{align*}
\rdF \varGamma(\etSite{X}, \mathscr{F}) \ttP[R] \rdF \varGamma(\etSite{Y}, \mathscr{G}) & \isoTo \rdF \varGamma \bigl(\etSite{Z}, f^{\ast}\mathscr{F} \ttP[R] g^{\ast}\mathscr{G}\bigr) \,, \\
\rdF \varGamma_!(\etSite{X}, \mathscr{F}) \ttP[R] \rdF \varGamma_!(\etSite{Y}, \mathscr{G}) & \isoTo \rdF \varGamma_!\bigl(\etSite{Z}, f^{\ast}\mathscr{F} \ttP[R] g^{\ast}\mathscr{G}\bigr) \,.
\end{align*}
Moreover there are two exact sequences of $R$-modules
\begin{align*}
0 & \to \bigoplus_{i+j=n} \upH{i}(\etSite{X}, \mathscr{F}) \otimes_R \upH{j}(\etSite{Y}, \mathscr{G}) \to \upH{n}\bigl(\etSite{Z}, f^{\ast}\mathscr{F} \ttP[R] g^{\ast}\mathscr{G}\bigr) \\
& \to \bigoplus_{i+j=n+1} \Tor^R_1\bigl(\upH{i}(\etSite{X}, \mathscr{F}), \upH{j}(\etSite{Y}, \mathscr{G})\bigr) \to 0 \,, \\
0 & \to \bigoplus_{i+j=n} \upHc{i}(\etSite{X}, \mathscr{F}) \otimes_R \upHc{j}(\etSite{Y}, \mathscr{G}) \to \upHc{n}\bigl(\etSite{Z}, f^{\ast}\mathscr{F} \ttP[R] g^{\ast}\mathscr{G}\bigr) \\
& \to \bigoplus_{i+j=n+1} \Tor^R_1\bigl(\upHc{i}(\etSite{X}, \mathscr{F}), \upHc{j}(\etSite{Y}, \mathscr{G})\bigr) \to 0 \,.
\end{align*}
\end{Thm}

\section{The functor \texorpdfstring{$\rdF f^!$}{f!} and The Trace Morphisms from SGA 4 \& \texorpdfstring{4$\frac{1}{2}$}{$4.5$}}\label{Sec5}

\begin{Prop}\label{195P2}
Let $f \colon X \to Y$ be a compactifiable morphism of schemes satisfying {\upshape(\dag)} such that all fibers of $f$ are of dimensions $\leqslant d$. Then for each $a \in \mathbb{Z}$, $\rdF f^!$ sends $\conD[\geqslant a]{\etSite{X}, R_{\smdot}}$ to $\conD[\geqslant a-2d]{\etSite{X}, R_{\smdot}}$.
\end{Prop}

\begin{proof}
See \cite[XVIII, 3.1.7]{SGA4}.
\end{proof}

\begin{Lem}
Let $f \colon X \to Y$ be a compactifiable morphism of schemes satisfying {\upshape(\dag)}. Then for every pair of objects $\mathscr{F}$ and $\mathscr{G}$ in $\conDn{\etSite{Y}, R_{\smdot}}$, there is a natural morphism
\[ \rdF f^! \mathscr{F} \ttP[R] f^{\ast}\mathscr{G} \to \rdF f^!(\mathscr{F} \ttP[R] \mathscr{G}) \]
in $\conDn{\etSite{X}, R_{\smdot}}$ which is functorial in $\mathscr{F}$ and $\mathscr{G}$.
\end{Lem}

\begin{Prop}\label{195P13}
Let $f \colon X \to Y$ and $g \colon Y \to Z$ two compactifiable morphisms of schemes satisfying {\upshape(\dag)}. For every pair of objects $\mathscr{F}$ and $\mathscr{G}$ in $\conDn{\etSite{Z}, R_{\smdot}}$, there is a natural morphism
\[ \rdF f^! \circ g^{\ast}\mathscr{F} \ttP[R] f^{\ast} \circ \rdF g^! \mathscr{G} \to \rdF (g \circ f)^!(\mathscr{F} \ttP[R] \mathscr{G}) \]
in $\conD{\etSite{X}, R_{\smdot}}$ which is functorial in $\mathscr{F}$ and $\mathscr{G}$.
\end{Prop}

\begin{proof}
We have
\begin{align*}
\rdF f^! \circ g^{\ast}\mathscr{F} \ttP[R] f^{\ast} \circ \rdF g^! \mathscr{G} & \to \rdF f^!\bigl(g^{\ast}\mathscr{F} \ttP[R] \rdF g^! \mathscr{G}\bigr) \\
& \to \rdF f^! \circ \rdF g^!(\mathscr{F} \ttP[R] \mathscr{G}) \\
& \isoTo \rdF (g \circ f)^!(\mathscr{F} \ttP[R] \mathscr{G}) \,. \qedhere
\end{align*}
\end{proof}

\begin{Prop}\label{195P15}
Let
\[ \xymatrix{X' \ar[r]^{f'} \ar[d]_{p} \CartS & Y' \ar[d]^{q} \\ X \ar[r]_{f} & Y} \]
be a Cartesian square of schemes satisfying {\upshape(\dag)}. Assume that $f$ is compactifiable.
\begin{enumerate}
\item For each object $\mathscr{F}$ in $\conD{\etSite{X}, R_{\smdot}}$, there is a natural morphism in $\conD{\etSite{Y'}, R_{\smdot}}$
\[ q^{\ast} \circ \rdF f_{\ast} \mathscr{F} \to \rdF f'_{\ast} \circ p^{\ast} \mathscr{F} \,. \]
\item For each object $\mathscr{G}$ in $\conD{\etSite{Y}, R_{\smdot}}$, there is a natural morphism in $\conD{\etSite{X'}, R_{\smdot}}$
\[ p^{\ast} \circ \rdF f^! \mathscr{G} \to \rdF f'^! \circ q^{\ast} \mathscr{G} \,. \]
\item Assume that $Y$ is an algebraic scheme over a field $k$, and there exists a $k$-scheme $T$ such that $Y' = Y \times_k T$. Then the morphisms in (1) and (2) are both isomorphisms.
\item For each object $\mathscr{G}$ in $\conD{\etSite{Y'}, R_{\smdot}}$, there is a natural morphism in $\conD{\etSite{X}, R_{\smdot}}$
\[ \rdF p_{\ast} \circ \rdF f'^! \mathscr{G} \isoTo \rdF f^! \circ \rdF q_{\ast} \mathscr{G} \,. \]
\end{enumerate}
\end{Prop}

\begin{proof}
(1) is induced by the classical base change morphisms.

(2) is from \cite[XVIII, 3.1.14.2]{SGA4}.

(3) is by \cite[Th. finitude, 1.9]{MR0463174}.

(4) is by \cite[XVIII, 3.1.12.3]{SGA4}.
\end{proof}

Now we review the trace morphisms.

\begin{Defi}
A morphism $f \colon X \to Y$ of schemes is said to be \emph{flat at dimension $d$} if there exists a nonempty open subset $U$ of $X$ satisfying the following conditions:
\begin{enumerate}
\item $f \colon U \to Y$ is flat;
\item for each point $y \in Y$, $U_y$ is either empty or $d$-dimensional;
\item every fiber of $X \setminus U \to Y$ is of dimension $< d$.
\end{enumerate}
\end{Defi}

By \cite[XVIII, 2.9]{SGA4}, for every compactifiable morphism $f \colon X \to Y$ of schemes satisfying (\dag) which is flat at dimension $d$, and for every object $\mathscr{G}$ in $\conD{\etSite{Y}, R_{\smdot}}$, we have a \emph{trace morphism}:
\[ \Tr_f \colon \rdF f_! \circ f^{\ast}\mathscr{G}\dangle{d} \to \mathscr{G} \,. \]
Since $\rdF f^!$ is right adjoint to $\rdF f_!$, the morphism $\Tr_f$ induces a canonical morphism in $\conD{\etSite{X}, R_{\smdot}}$:
\[ \mathrm{t}_f \colon f^{\ast}\mathscr{G}\dangle{d} \to \rdF f^!\mathscr{G} \,. \]
Moreover we have a commutative diagram
\begin{equation}\label{Eq:4c}
\vcenter{\xymatrix@C+2em{\rdF f_! \circ f^{\ast}\mathscr{G}\dangle{d} \ar[r]^-{\rdF f_!(\mathrm{t}_f)} \ar[rd]_-{\Tr_f} & \rdF f_! \circ \rdF f^! \mathscr{G} \ar[d] \\ & \mathscr{G}}}
\end{equation}

By \cite[XVIII, 3.2.5]{SGA4}, we have

\begin{Prop}\label{195T5}
Let $f \colon X \to Y$ be a compactifiable smooth morphism of relative dimension $d$ of schemes satisfying {\upshape(\dag)}. Then for any object $\mathscr{G}$ in $\conD{\etSite{Y}, R_{\smdot}}$, the canonical morphism
\[ \mathrm{t}_f \colon f^{\ast}\mathscr{G}\dangle{d} \isoTo \rdF f^!\mathscr{G} \]
is an isomorphism in $\conD{\etSite{X}, R_{\smdot}}$.
\end{Prop}

The following propositions \ref{195P11}-\ref{195p4a} are deduced from \cite[XVIII, 2.9]{SGA4}.

\begin{Prop}\label{195P11}
Let
\[ \xymatrix{X' \ar[r]^{f'} \ar[d]_{p} \CartS & Y' \ar[d]^{q} \\ X \ar[r]_{f} & Y} \]
be a Cartesian square of schemes satisfying {\upshape(\dag)}. Assume that $f$ is compactifiable and flat at dimension $d$. Then $f'$ is also flat at dimension $d$, and for each object $\mathscr{G}$ in $\conD{\etSite{Y}, R_{\smdot}}$ we have
\begin{enumerate}
\item the composite morphism
\[ (\rdF f'_!) \circ f'^{\ast} \circ q^{\ast} \mathscr{G}\dangle{d} = (\rdF f'_!) \circ p^{\ast} \circ f^{\ast} \mathscr{G}\dangle{d} \isoFrom q^{\ast} \circ (\rdF f_!) \circ f^{\ast} \mathscr{G}\dangle{d} \xrightarrow{q^{\ast}(\Tr_f)} q^{\ast}\mathscr{G} \]
is equal to $\Tr_{f'}$;
\item the composite morphism
\[ f'^{\ast} \circ q^{\ast}\mathscr{G}\dangle{d} = p^{\ast} \circ f^{\ast}\mathscr{G}\dangle{d} \xrightarrow{p^{\ast}(\mathrm{t}_f)} p^{\ast} \circ \rdF f^! \mathscr{G} \to \rdF f'^! \circ q^{\ast} \mathscr{G} \]
is equal to $\mathrm{t}_{f'}$, where the last morphism is defined in Proposition \ref{195P15} (2).
\end{enumerate}
\end{Prop}

\begin{Prop}\label{195P4}
Let $f \colon X \to Y$ and $g \colon Y \to Z$ two compactifiable morphisms of schemes satisfying {\upshape(\dag)} which are flat at dimension $d$ and $e$ respectively, $\mathscr{H}$ an object in $\conD{\etSite{Z}, R_{\smdot}}$. Then we have
\begin{enumerate}
\item The composite morphism
\begin{gather*}
\rdF (g \circ f)_! \circ (g \circ f)^{\ast}\mathscr{H}\dangle{d+e} \isoTo (\rdF g_!) \circ (\rdF f_!) \circ f^{\ast} \circ g^{\ast} \mathscr{H}\dangle{d+e} \\
\xrightarrow{\rdF g_!(\Tr_f)} (\rdF g_!) \circ g^{\ast} \mathscr{H}\dangle{e} \xrightarrow{\Tr_g} \mathscr{H}
\end{gather*}
is equal to $\Tr_{g \circ f}$.
\item The composite morphism
\[ f^{\ast}g^{\ast}\mathscr{H}\dangle{d+e} \xrightarrow{\mathrm{t}_f} \rdF f^! \circ g^{\ast}\mathscr{H}\dangle{e} \xrightarrow{\rdF f^!(\mathrm{t}_g)} \rdF f^! \circ \rdF g^! \mathscr{H} \isoTo \rdF (g \circ f)^! \mathscr{H} \]
is equal to $\mathrm{t}_{g \circ f}$.
\end{enumerate}
\end{Prop}

\begin{Prop}\label{195p4a}
Let $f \colon X \to Y$ be a finite morphism of schemes satisfying {\upshape(\dag)} such that $f_{\ast}\mathcal{O}_X$ is a locally free $\mathcal{O}_Y$-module of degree $d$. Then for each object $\mathscr{F}$ in $\conD{\etSite{Y}, R_{\smdot}}$, the composite morphism
\[ \mathscr{F} \to f_{\ast}f^{\ast}\mathscr{F} \xrightarrow{\Tr_f} \mathscr{F} \]
is equal to the multiplication by $n$.
\end{Prop}

The following proposition show that the trace morphism is essentially determinated by the generic points. Let $A$ be a Noetherian ring (in particular $A = R_n$).

\begin{Prop}\label{3hj01r}
Let $X$ be a $n$-dimensional algebraic scheme over $k$, $X_1,X_2,\ldots,X_r$ all irreducible components of dimension $n$ of $X$, $F$ an $A$-module. For each $i$, let $Y_i \neq \varnothing$ be an open subset of $X$ contained $X_i \setminus \bigcup\limits_{j \neq i}X_j$ and regard $Y_i$ as a reduced subscheme of $X$. For each $i$, let $x_i$ be the generic point of $X_i$ and put $a_i \coloneqq \mathrm{length}(\mathcal{O}_{X,x_i})$. Then there is a canonical isomorphism $\omega$ of $A$-modules which makes a commutative diagram.
\[ \xymatrix{\bigoplus\limits^r_{i=1} \upHc{2n}\bigl(Y_{i,\textup{\'et}}, F(n)\bigr) \ar[rr]^-{\omega}_-{\cong} \ar[rd]_{\bigoplus\limits^r_{i=1} a_i \cdot \Tr_i} & & \upHc{2n}\bigl(\etSite{X}, F(n)\bigr) \ar[ld]^{\Tr} \\ & F} \]
\end{Prop}

\section{\texorpdfstring{$\ell$}{l}-adic homology for morphisms of algebraic schemes}\label{Sec6}

\begin{Nota}
Let $f \colon X \to Y$ be a compactifiable morphism of schemes satisfying (\dag). For each object $\mathscr{N}$ in $\conD{\etSite{Y}, R_{\smdot}}$ and for each $n \in \mathbb{Z}$, we define the \emph{$n$-th $\ell$-adic homology associated to $f$} to be
\begin{align*}
	\rmH{n}(X \xrightarrow{f} Y, \mathscr{N})  & \coloneqq \upH{-n}(\etSite{X}, \, \rdF f^! \mathscr{N}) \\
	&= \Hom_{\conD{\etSite{X}, R_{\smdot}}}\bigl(R_X, \, \rdF f^! \mathscr{N}[-n]\bigr) \,,
\end{align*}
which is an $R$-module.

For convenient to define pull-backs along flat morphisms and cycle maps, we also define
\begin{align*}
\mathbb{H}_n(X \xrightarrow{f} Y, \mathscr{N}) & \coloneqq \rmH{2n}\bigl(X \xrightarrow{f} Y, \mathscr{N}(-n)\bigr) \\
&= \Hom_{\conD{\etSite{X}, R_{\smdot}}}\bigl(R_X, \, \rdF f^! \mathscr{N}\dangle{-n}\bigr) \,. \\
\end{align*}

We set
\begin{align*}
	\rmH{\ast}(X \xrightarrow{f} Y, \mathscr{N}) & \coloneqq \bigoplus_{n \in \mathbb{Z}}\rmH{n}(X \xrightarrow{f} Y, \mathscr{N}) \,, \\
	\mathbb{H}_{\ast}(X \xrightarrow{f} Y, \mathscr{N}) & \coloneqq \bigoplus_{n \in \mathbb{Z}}\mathbb{H}_n(X \xrightarrow{f} Y, \mathscr{N}) \,.
\end{align*}

We also use $\rmH{n}(X/Y, \mathscr{N})$ to denote $\rmH{\ast}(X \xrightarrow{f} Y, \mathscr{N})$ if no confusion arise. Similarly we may define $\mathbb{H}_n(X/Y, \mathscr{N})$.

\smallskip

If $X$ is an algebraic schemes over a separably closed field $k$ and $N$ is an object in $\fgD{R}$, we write
\[ \rmH{n}(X, N) \coloneqq \rmH{n}(X \to \Spec k, N) \,, \qquad \mathbb{H}_n(X, N) \coloneqq \mathbb{H}_n(X \to \Spec k, N) \,. \]
\end{Nota}

By Proposition \ref{195T5} we have

\begin{Lem}\label{195T5t}
Let $f \colon X \to Y$ be a compactifiable smooth morphism of relative dimension $d$ of schemes satisfying {\upshape(\dag)}. Then for each object $\mathscr{N}$ in $\conD{\etSite{Y}, R_{\smdot}}$ and for $n \in \mathbb{Z}$, the morphism $\mathrm{t}_f$ induces a canonical isomorphism of $R$-modules:
\[ \mathbb{H}^{d-n}(X,  f^{\ast}\!\mathscr{N}) \isoTo \mathbb{H}_n(X \xrightarrow{f} Y, \mathscr{N}) \,. \]
\end{Lem}

\begin{Prop}\label{T-P0}
	Let $f \colon X \to S$ be a compactifiable morphism of schemes satisfying {\upshape(\dag)}, $Y$ a closed subscheme of $X$ and $U \coloneqq X \setminus Y$. Then we have a long exact sequence
	\[ \cdots \to \rmH{n}(Y/S, \mathscr{N}) \to \rmH{n}(X/S, \mathscr{N}) \to \rmH{n}(U/S, \mathscr{N}) \to \rmH{n-1}(Y/S, \mathscr{N}) \to \cdots \,. \]
\end{Prop}

\begin{proof}
	Put $\mathscr{M} \coloneqq \rdF f^!\mathscr{N}$. Then the proposition follows from the distinguished triangle
	\[ i_{\ast}i^!\mathscr{M} \to \mathscr{M} \to j_{\ast}j^{\ast}\mathscr{M} \to i_{\ast}i^!\mathscr{M}[1] \,, \]
	where $i \colon Y \hookrightarrow X$ and $j \colon U \hookrightarrow X$ are the inclusions.
\end{proof}

\begin{Prop}[Mayer-Vietoris Sequence]\label{T-P1}
	Let $f \colon X \to S$ be a compactifiable morphism of schemes satisfying {\upshape(\dag)}, $X_1$ and $X_2$ two closed subschemes of $X$ such that $X = X_1 \cup X_2$ (as sets). Then we have a long exact sequence
	\begin{align*}
		\cdots & \to \rmH{n}\bigl((X_1 \cap X_2)/S, \mathscr{N}\bigr) \to \rmH{n}(X_1/S, \mathscr{N}) \oplus \rmH{n}(X_2/S, \mathscr{N}) \to \rmH{n}(X/S, \mathscr{N}) \\
		& \to \rmH{n-1}\bigl((X_1 \cap X_2)/S, \mathscr{N}\bigr) \to \cdots \,.
	\end{align*}
\end{Prop}

\begin{proof}
	Put $\mathscr{M} \coloneqq \rdF f^!\mathscr{N}$. Then the proposition follows from the distinguished triangle
	\[ i_{\ast} \circ \rdF i^! \mathscr{M} \to i_{1,\ast} \circ \rdF i_1^! \mathscr{M} \oplus i_{2,\ast} \circ \rdF i_2^! \mathscr{M} \to \mathscr{M} \to i_{\ast} \circ \rdF i^! \mathscr{M}[1] \,, \]
	where $i \colon X_1 \cap X_2 \hookrightarrow X$, $i_1 \colon X_1 \hookrightarrow X$, $i_2 \colon X_2 \hookrightarrow X$ are the inclusions.
\end{proof}

\begin{Prop}[Vanishing]\label{195P9}
Let $f \colon X \to Y$ be a compactifiable morphism of schemes satisfying {\upshape(\dag)} such that all fibers of $f$ are of dimensions $\leqslant d$, $\mathscr{N}$ an object in $\conD[\geqslant a]{\etSite{Y}, R_{\smdot}}$. Then $\rmH{n}(X/Y, \mathscr{N}) = 0$ whenever $n > 2d - a$.
\end{Prop}

\begin{proof}
By Proposition \ref{195P2}, $\rdF f^!\mathscr{N}[-n] \in \conD[\geqslant a-2d+n]{\etSite{X}, R_{\smdot}}$. Thus if $a - 2d + n > 0$, then
\[ \rmH{n}(X \xrightarrow{f} Y, \mathscr{N}) = \Hom_{\conD{\etSite{X}, R_{\smdot}}}\bigl(R_X, \rdF f^! \mathscr{N}[-n]\bigr) = 0 \,. \qedhere \]
\end{proof}

\begin{Prop}\label{78fu0i}
Let $f \colon X \to S$ be a compactifiable morphism of schemes satisfying {\upshape(\dag)}, $Y$ a closed subscheme of $X$  such that $\dim Y_s \leqslant d$ for all $s \in S$, $X' \coloneqq X \setminus Y$, $\mathscr{N}$ an object in $\conD[\geqslant a]{\etSite{S}, R_{\smdot}}$. Then for each integer $n > 2d + 1 - a$, there is a canonical isomorphism of $R$-modules
\[ \rmH{n}(X/S, \mathscr{N}) \isoTo \rmH{n}(X'/S, \mathscr{N}) \,. \]
\end{Prop}

\begin{proof}
	Apply Proposition \ref{T-P0} and Proposition \ref{195P9}.
\end{proof}

\begin{Nota}
Let $f \colon X \to Y$ be a compactifiable morphism of schemes satisfying (\dag). For each object $\mathscr{G}$ in $\conD{\etSite{Y}, R_{\smdot}}$, we define
\[ \delta_f \colon \mathscr{G} \to \rdF f_{\ast} \circ f^{\ast}\mathscr{G} \quad \text{and} \quad \theta_f \colon \rdF f_! \circ \rdF f^! \mathscr{G} \to \mathscr{G} \]
to be the canonical morphisms induced by the adjunctions $f^{\ast} \dashv \rdF f_{\ast}$ and $\rdF f_! \mapsto \rdF f^!$ respectively.
\end{Nota}

The following map is a kind of variant of the Gysin homomorphism.

\begin{Defi}[Push-forward]
Let $p \colon X \to S$ and $q \colon Y \to S$ be two compactifiable morphisms of schemes satisfying {\upshape(\dag)}, $f \colon X \to Y$ a proper $S$-morphism. For every object $\mathscr{N}$ in $\conD{\etSite{S}, R_{\smdot}}$ and for every $n \in \mathbb{Z}$, we define a homomorphism of $R$-modules
\[ f_{\ast} \colon \rmH{n}(X/S, \mathscr{N}) \to \rmH{n}(Y/S, \mathscr{N}) \]
as follows. For each $\alpha \in \rmH{n}(X/S, \mathscr{N})$, $f_{\ast}(\alpha)$ is defined to be the composition
\[ R_Y \xrightarrow{\delta_f} \rdF f_{\ast} R_X \xrightarrow{\rdF f_{\ast}(\alpha)} \rdF f_{\ast} \circ \rdF p^! \mathscr{N}[-n] \isoTo \rdF f_{\ast} \circ \rdF f^! \circ \rdF q^! \mathscr{N}[-n] \xrightarrow{\theta_f} \rdF q^! \mathscr{N}[-n] \,. \]
\end{Defi}

\begin{Prop}\label{195P8}
Let $X \xrightarrow{f} Y \xrightarrow{g} Z \xrightarrow{h} S$ be a sequence of morphisms of schemes satisfying {\upshape(\dag)} such that $f$ and $g$ are proper, and $h$ is compactifiable. Then for all $\mathscr{N} \in \conD{\etSite{S}, R_{\smdot}}$ and $n \in \mathbb{Z}$, we have
\[ (g \circ f)_{\ast} = g_{\ast} \circ f_{\ast} \colon \rmH{n}(X/S, \mathscr{N}) \to \rmH{n}(Z/S, \mathscr{N}) \,. \]
\end{Prop}

\begin{proof}
This is by the following simple lemma.
\end{proof}

\begin{Lem}
Let $f \colon X \to Y$ and $g \colon Y \to Z$ be two compactifiable morphisms of schemes satisfying {\upshape(\dag)}, $\mathscr{H}$ an object in $\conD{\etSite{Z}, R_{\smdot}}$. Then we have
\begin{enumerate}
\item The following composition is equal to $\delta_{g \circ f}$
\[ \mathscr{H} \xrightarrow{\delta_g} \rdF g_{\ast} \circ g^{\ast}\mathscr{H} \xrightarrow{\rdF g_{\ast}(\delta_f)} \rdF g_{\ast} \circ \rdF f_{\ast} \circ f^{\ast} \circ g^{\ast}\mathscr{H} \isoTo \rdF (g \circ f)_{\ast} \circ (g \circ f)^{\ast}\mathscr{H} \,. \]
\item The following composition is equal to $\theta_{g \circ f}$
\[ \rdF (g \circ f)_! \circ \rdF (g \circ f)^! \mathscr{H} \isoTo \rdF g_! \circ \rdF f_! \circ \rdF f^! \circ \rdF g^! \mathscr{H} \xrightarrow{\rdF g_!(\theta_f)} \rdF g_! \circ \rdF g^! \mathscr{H} \xrightarrow{\theta_g} \mathscr{H} \,. \]
\end{enumerate}
\end{Lem}

\begin{Defi}[Pull-back]
Let $p \colon X \to S$ and $q \colon Y \to S$ be two compactifiable morphisms of schemes satisfying {\upshape(\dag)}, $f \colon X \to Y$ an $S$-morphism which is flat at dimension $d$. For every object $\mathscr{N}$ in $\conD{\etSite{S}, R_{\smdot}}$ and for every $n \in \mathbb{Z}$, we define a homomorphism of $R$-modules
\[ f^{\ast} \colon \mathbb{H}_n(Y/S, \mathscr{N}) \to \mathbb{H}_{n+d}(X/S, \mathscr{N}) \]
as follows. For each $\beta \in \mathbb{H}_n(Y/S, \mathscr{N})$, $f^{\ast}(\beta)$ is defined to be the composition
\[ R_X \xrightarrow{\mathrm{t}_f} \rdF f^! R_Y\dangle{-d} \xrightarrow{\rdF f^!(\beta)} \rdF f^! \circ \rdF q^! \mathscr{N}\dangle{-(n+d)} \isoTo \rdF p^!\mathscr{N}\dangle{-(n+d)} \,. \]
\end{Defi}

\begin{Prop}\label{217P4}
Let $X \xrightarrow{f} Y \xrightarrow{g} Z \to S$ be a sequence of compactifiable morphisms of schemes satisfying {\upshape(\dag)} such that $f$ and $g$ are flat at dimension $d$ and $e$ respectively. Then for all $\mathscr{N} \in \conD{\etSite{S}, R_{\smdot}}$ and $n \in \mathbb{Z}$, we have
\[ (g \circ f)^{\ast} = f^{\ast} \circ g^{\ast} \colon \mathbb{H}_n(X/S, \mathscr{N}) \to \mathbb{H}_{n+d+e}(Z/S, \mathscr{N}) \,. \]
\end{Prop}

\begin{proof}
This follows from Proposition \ref{195P4} (2).
\end{proof}

\begin{Thm}\label{217P3}
Let $S$ be a scheme satisfying {\upshape(\dag)}, $r \colon Y \to S$ a compactifiable morphism. Let
\[ \xymatrix{X' \ar[r]^{f'} \ar[d]_{p} \CartS & Y' \ar[d]^{q} \\ X \ar[r]_{f} & Y} \]
be a Cartesian square of schemes such that $f$ is proper and $q$ is compactifiable and flat at dimension $d$, $\mathscr{N}$ an object in $\conD{\etSite{S}, R_{\smdot}}$ and $n \in \mathbb{Z}$. Then we have
\[ q^{\ast} \circ f_{\ast} = f'_{\ast} \circ p^{\ast} \colon \mathbb{H}_n(X/S, \mathscr{N}) \to \mathbb{H}_{n+d}(Y'/S, \mathscr{N}) \,. \]
\end{Thm}

\begin{proof}
Put $\mathscr{M} \coloneqq \rdF r^! \mathscr{N}$. Let $\alpha \in \mathbb{H}_n(X/S, \mathscr{N})$. Then $q^{\ast} \circ f_{\ast}(\alpha)$ is equal to the composition
\begin{gather*}
R_{Y'} \xrightarrow{\mathrm{t}_q} \rdF q^! R_Y\dangle{-d} \xrightarrow{\rdF q^!(\delta_f)} \rdF q^! \circ \rdF f_{\ast} R_X\dangle{-d} \xrightarrow{\rdF q^! \circ \rdF f_{\ast}(\alpha)} \\
\rdF q^! \circ \rdF f_{\ast} \circ \rdF f^! \mathscr{M}\dangle{-(n+d)} \xrightarrow{\rdF q^!(\theta_f)} \rdF q^!\mathscr{M}\dangle{-(n+d)} \,;
\end{gather*}
and $f'_{\ast} \circ p^{\ast}(\alpha)$ is equal to the composition
\begin{gather*}
R_{Y'} \xrightarrow{\delta_{f'}} \rdF f'_{\ast} \circ f'^{\ast} R_{Y'} = \rdF f'_{\ast} \circ p^{\ast} R_X \xrightarrow{\rdF f'_{\ast}(\mathrm{t}_p)} \rdF f'_{\ast} \circ \rdF p^! R_{X}\dangle{-d} \xrightarrow{\rdF f'_{\ast} \circ \rdF p^!(\alpha)} \\
\rdF f'_{\ast} \circ \rdF p^! \circ \rdF f^!\mathscr{M}\dangle{-(n+d)} = \rdF f'_{\ast} \circ \rdF f'^! \circ \rdF q^! \mathscr{M}\dangle{-(n+d)} \xrightarrow{\theta_{f'}} \rdF q^! \mathscr{M}\dangle{-(n+d)} \,.
\end{gather*}
After applying Proposition \ref{195P11} (2) to $\mathrm{t}_p$, we obtain that the morphism $f'_{\ast} \circ p^{\ast}(\alpha)$ is equal to the composition
\begin{gather*}
R_{Y'} \xrightarrow{\delta_{f'}} \rdF f'_{\ast} \circ f'^{\ast} R_{Y'} = \rdF f'_{\ast} \circ f'^{\ast} \circ q^{\ast} R_Y \xrightarrow{\rdF f'_{\ast} \circ f'^{\ast}(\mathrm{t}_q)} \rdF f'_{\ast} \circ f'^{\ast} \circ q^! R_Y\dangle{-d} \\
\to \rdF f'_{\ast} \circ \rdF p^! \circ f^{\ast} R_Y\dangle{-d} = \rdF f'_{\ast} \circ \rdF p^! R_X\dangle{-d} \xrightarrow{\rdF f'_{\ast} \circ \rdF p^!(\alpha)} \\
\rdF f'_{\ast} \circ \rdF p^! \circ \rdF f^!\mathscr{M}\dangle{-(n+d)} = \rdF f'_{\ast} \circ \rdF f'^! \circ \rdF q^! \mathscr{M}\dangle{-(n+d)} \xrightarrow{\theta_{f'}} \rdF q^! \mathscr{M}\dangle{-(n+d)} \,.
\end{gather*}
Consider the following diagram
\begin{gather*}
\hspace{1.2em}\xymatrix@C+2em{R_{Y'} \ar[r]^-{\mathrm{t}_q} \ar[d]_{\delta_{f'}} \ar@{}[dr]|-*+{\circlearrowright} & \rdF q^! R_Y\dangle{-d} \ar[d]|{\delta_{f'}} \ar[r]^-{\rdF q^!(\delta_f)} \ar@{}[dr]|-*+\txt{(a)} & *+!<2.4em,0ex>{\rdF q^! \circ \rdF f_{\ast} R_X\dangle{-d}} \ar@<-2em>[d]|{\cong} \ar[r]^-{\rdF q^! \circ \rdF f_{\ast}(\alpha)} \ar@{}[dr]|-*+{\circlearrowright} & *+!<3em,0ex>{\ } \\
  \rdF f'_{\ast} \circ f'^{\ast} R_{Y'} \ar[r]_-{\rdF f'_{\ast} \circ f'^{\ast}(\mathrm{t}_q)} & \rdF f'_{\ast} \circ f'^{\ast} \circ \rdF q^! R_Y\dangle{-d} \ar[r] & *+!<3em,0ex>{\rdF f'_{\ast} \circ \rdF p^! \circ f^{\ast} R_Y\dangle{-d}} \ar[r]_-{\rdF f'_{\ast} \circ \rdF p^!(\alpha)} & *+!<3em,0ex>{\ }} \\
\xymatrix{\rdF q^! \circ \rdF f_{\ast} \circ \rdF f^! \mathscr{M}\dangle{-(n+d)} \ar[d]|{\cong} \ar[rr]^-{\rdF q^!(\theta_f)} & \ar@{}[d]|-*+\txt{(b)} & *+!<2.3em,0ex>{\rdF q^! \mathscr{M}\dangle{-(n+d)}} \ar@{=}@<-2.4em>[d] \\
  \rdF f'_{\ast} \circ \rdF p^! \circ \rdF f^! \mathscr{M}\dangle{-(n+d)} \ar@{=}[r] & *+!<1.5em,0ex>{\rdF f'_{\ast} \circ \rdF f'^! \circ \rdF q^! \mathscr{M}\dangle{-(n+d)}} \ar[r]_-{\theta_{f'}} & *+!<2.3em,0ex>{\rdF q^! \mathscr{M}\dangle{-(n+d)}}}
\end{gather*}
\newcommand{\inMatrix}[1]{\begingroup%
  \def\objectstyle{\scriptscriptstyle}%
  $\vcenter{\xymatrix@C=1em@R=1.8ex{#1}}$%
\endgroup}%
where $\circlearrowright$ means commutative square. The commutativity of (a) and (b) are by the following simple Lemma \ref{195L4}. So the whole diagram is commutative. Note that the composition along the direction \inMatrix{\bullet \ar[r] & \bullet \ar[d] \\ & \bullet} in above diagram is equal to $q^{\ast} \circ f_{\ast}(\alpha)$; and the composition along \inMatrix{\bullet \ar[d] \\ \bullet \ar[r] & \bullet} is equal to $f'_{\ast} \circ p^{\ast}(\alpha)$. Thus $q^{\ast} \circ f_{\ast}(\alpha) = f'_{\ast} \circ p^{\ast}(\alpha)$.
\end{proof}

\begin{Lem}\label{195L4}
Let
\[ \xymatrix{X' \ar[r]^{f'} \ar[d]_{p} \CartS & Y' \ar[d]^{q} \\ X \ar[r]_{f} & Y} \]
be a Cartesian square of schemes satisfying {\upshape(\dag)} with all morphisms compactifiable. Then we have
\begin{enumerate}
\item For each object $\mathscr{G}$ in $\conD{\etSite{Y}, R_{\smdot}}$, the diagram
\[ \xymatrix@C+2em{\rdF q^! \mathscr{G} \ar[r]^-{\rdF q^!(\delta_f)} \ar[d]_{\delta_{f'}} & \rdF q^! \circ \rdF f_{\ast} \circ f^{\ast}\mathscr{G} \ar[d]^{\varphi}_{\cong} \\
  \rdF f'_{\ast} \circ f'^{\ast} \circ \rdF q^!\mathscr{G} \ar[r]^-{\rdF f'_{\ast}(\psi)} & \rdF f'_{\ast} \circ \rdF p^! \circ f^{\ast} \mathscr{G}} \]
is commutative in $\conD{\etSite{X'}, R_{\smdot}}$, where $\varphi$ is defined in Proposition \ref{195P15} (4) and $\psi$ is defined in Proposition \ref{195P15} (2).
\item Assume that $f$ is proper. Then for each object $\mathscr{G}$ in $\conD{\etSite{Y}, R_{\smdot}}$, the diagram
\[ \xymatrix{\rdF f'_{\ast} \circ \rdF p^! \circ \rdF f^! \mathscr{G} \ar[d]_{\rdF f'_{\ast}(\beta)}^{\cong} \ar[r]^-{\alpha}_-{\cong} & \rdF q^! \circ \rdF f_{\ast} \circ \rdF f^! \mathscr{G} \ar[d]^{\rdF q^!(\theta_f)} \\ \rdF f'_{\ast} \circ \rdF f'^! \circ \rdF q^! \mathscr{G} \ar[r]_-{\theta_{f'}} & \rdF q^! \mathscr{G}} \]
is commutative in $\conD{\etSite{Y'}, R_{\smdot}}$, where $\alpha$ is defined in Proposition \ref{195P15} (4) and $\beta$ is induced by the composition
\[ \rdF p^! \circ \rdF f^! \isoTo \rdF (f \circ p)^! = \rdF (q \circ f')^! \isoTo \rdF f'^! \circ \rdF q^! \,. \]
\end{enumerate}
\end{Lem}

\begin{Defi}[Base Change]
Let
\[ \xymatrix{X' \ar[r]^{f'} \ar[d]_{p} \CartS & S' \ar[d]^{u} \\ X \ar[r]_{f} & S} \]
be a Cartesian square of schemes satisfying (\dag) with $f$ compactifiable. For every object $\mathscr{N}$ in $\conD{\etSite{S}, R_{\smdot}}$ and for every $n \in \mathbb{Z}$, we define a homomorphism of $R$-modules
\[ u^{\ast} \colon \rmH{n}(X/S, \mathscr{N}) \to \rmH{n}(X'/S', u^{\ast}\mathscr{N}) \]
as follows. For each $\alpha \in \rmH{n}(X/S, \mathscr{N})$, $u^{\ast}(\alpha)$ is defined to be the composition
\[ R_{X'} = p^{\ast}R_X \xrightarrow{p^{\ast}(\alpha)} p^{\ast} \circ \rdF f^! \mathscr{N}[-n] \xrightarrow{\varphi} \rdF f'^! \circ u^{\ast} \mathscr{N}[-n] \,, \]
where $\varphi$ is defined in Proposition \ref{195P15} (2).
\end{Defi}

We have the following three obvious propositions about the base change homomorphisms.

\begin{Prop}\label{base0}
Let $k \subseteq K$ be two separably closed fields, $f \colon X \to S$ a morphism of algebraic schemes over $k$, $u \colon S_K \to S$ the projection. Then for each object $\mathscr{N}$ in $\conD{\etSite{S}, R_{\smdot}}$ and for each $n \in \mathbb{Z}$, the homomorphism
\[ u^{\ast} \colon \rmH{n}(X/S, \mathscr{N}) \isoTo \rmH{n}(X_K/S_K, u^{\ast}\mathscr{N}) \]
is an isomorphism.
\end{Prop}

\begin{proof}
It follows from Proposition \ref{195P15} (3) and Theorem \ref{194P14}.
\end{proof}

\begin{Prop}
Let
\[ \xymatrix{X'' \ar[r] \ar[d] \CartS & X' \ar[r] \ar[d] \CartS & X \ar[d] \\ S'' \ar[r]_{v} & S' \ar[r]_{u} & S} \]
be a commutative diagram of schemes satisfying {\upshape(\dag)} with both squares Cartesian, and all three vertical arrows being compactifiable. Then for all $\mathscr{N} \in \conD{\etSite{S}, R_{\smdot}}$ and $n \in \mathbb{Z}$, we have
\[ (u \circ v)^{\ast} = v^{\ast} \circ u^{\ast} \colon \rmH{n}(X/S, \mathscr{N}) \to \rmH{n}(X''/S'', (u \circ v)^{\ast}\mathscr{N}) \,. \]
\end{Prop}

\begin{Prop}\label{j4i0yn}
Let
\[ \xymatrix{X' \ar[r]^{f'} \ar[d] \CartS & Y' \ar[r] \ar[d] \CartS & S' \ar[d]^{u} \\ X \ar[r]_{f} & Y \ar[r] & S} \]
be a commutative diagram of schemes satisfying {\upshape(\dag)} with both squares Cartesian, and all level arrows being compactifiable. Let $\mathscr{N}$ an object in $\conD{\etSite{S}, R_{\smdot}}$ and $n \in \mathbb{Z}$. Then we have
\begin{enumerate}
\item If $f$ is proper, then
\[ u^{\ast} \circ f_{\ast} = f'_{\ast} \circ u^{\ast} \colon \rmH{n}(X/S, \mathscr{N}) \to \rmH{n}(Y'/S', u^{\ast}\mathscr{N}) \,. \]
\item If $f$ is flat at dimension $d$, then
\[ u^{\ast} \circ f^{\ast} = f'^{\ast} \circ u^{\ast} \colon \mathbb{H}_n(Y/S, \mathscr{N}) \to \mathbb{H}_{n+d}(X'/S', u^{\ast}\mathscr{N}) \,. \]
\end{enumerate}
\end{Prop}

\begin{Defi}[Galois action]
Let $k_0$ be a field, $k$ the separably closed field of $k_0$, $G \coloneqq \mathrm{Gal}(k/k_0)$, $X$ an algebraic scheme over $k$, $Y_0$ an algebraic scheme over $k_0$, $Y \coloneqq Y_0 \otimes_{k_0} k$, $\mathscr{N}_0$ an object in $\conD{Y_{0,\textup{\'et}}, R_{\smdot}}$ and $\mathscr{N}$ the pull-back of $\mathscr{N}_0$ on $Y$. Then there is an action of $G$ on $\rmH{n}(X/Y, \mathscr{N})$ defined by
\[ (g, \alpha) \mapsto (\iDe_{Y_0} \otimes g)^{\ast}(\alpha) \,, \qquad g \in G, \, \alpha \in \rmH{n}(X/Y, \mathscr{N}) \,. \]

In paritcular if $N \in \fgD{R}$ and $n \in \mathbb{Z}$, then there is a Galois action of $G$ on $\rmH{n}(X, N)$.
\end{Defi}

The following theorem is used to prove that cycle maps eliminate algebraic equivalent classes.

\begin{Thm}\label{eaec8}
Let $f \colon X \to Y$ be a morphism of algebraic schemes over a separably closed field $k$, $Z$ a nonsingular variety over $k$, $\mathscr{N}$ an object in $\conD{\etSite{Y}, R_{\smdot}}$,
\[ \alpha \in \rmH{n}\bigl((X \times_k Z)/(Y \times_k Z), \mathrm{pr}_1^{\ast}\mathscr{N}\bigr) \,. \]
For each $z \in Z(k)$, put
\[ j_z \coloneqq \iDe_Y \times z \colon Y \to Y \times_k Z \,. \]
Then $z \mapsto j_z^{\ast}(\alpha)$ is a constant map from $Z(k)$ to $\rmH{n}(X/Y, \mathscr{N})$.
\end{Thm}

\begin{proof}
By Proposition \ref{base0} we may assume that $k$ is algebraically closed. Since every two rational points of $Z$ can be jointed by a series of nonsingular curves, we may further assume that $Z$ is a complete nonsingular curve. First we have a commutative diagram with both squares Cartesian.
\[ \xymatrix{X \times_k Z \ar[r]^-{f'} \ar[d]_{p} \CartS & Y \times_k Z \ar[r]^-{v} \ar[d]|{\vphantom{f}u} \CartS & Z \ar[d]^{r} \\
X \ar[r]_{f} & Y \ar[r]_{q} & \Spec k} \]
By Proposition \ref{195P15}, we have
\[ \rdF f'^! \circ u^{\ast}\mathscr{N}[-n] \isoTo p^{\ast} \rdF f^! \mathscr{N}[-n] = p^{\ast}\rdF f^! \mathscr{N}[-n] \ttP[R]  (v \circ f')^{\ast}R_Z \,. \]
 Since $Z$ is a complete nonsingular curve over $k$, we have $\upH{0}(\etSite{Z}, R) \cong R$, $\upH{1}(\etSite{Z}, R) \cong R^{\oplus g}$ (where $g$ is the genus of $Z$), and $\upH{2}(\etSite{Z}, R) \isoTo[\Tr] R$ are all free $R$-modules.  Now we apply Theorem \ref{194T2} to obtain an isomorphism:
\begin{gather*}
\rmH{n}\bigl((X \times_k Z)/(Y \times_k Z), \mathrm{pr}_1^{\ast}\mathscr{N}\bigr) \isoTo \upH{0}(X \times_k Z, \, p^{\ast}\rdF f^! \mathscr{N}[-n] \ttP[R]  (v \circ f')^{\ast}R_Z) \\
\isoTo \rmH{n}(X/Y,\mathscr{N}) \,{\textstyle\bigoplus} \bigl(\rmH{n+1}(X/Y,\mathscr{N}) \otimes_R \upH{1}(\etSite{Z}, R)\bigr) \,{\textstyle\bigoplus}\, \bigl(\rmH{n+2}(X/Y,\mathscr{N}) \otimes_R \upH{2}(\etSite{Z}, R)\bigr) \,.
\end{gather*}
Let $\beta \in \rmH{n}(X/Y,\mathscr{N})$ be the image of $\alpha$ induced by above isomorphism. Then $j_z^{\ast}(\alpha) = \beta$ for all $z \in Z(k)$.
\end{proof}

\section{The Cycle Maps for Chow Groups}\label{Sec7}

In this section, we construct the cycle maps for arbitrary algebraic schemes over $k$, where $k$ is a field such that $\Char k \neq \ell$ and $\mathrm{cd}_{\ell}(k) < \infty$.

\begin{Nota}
Let $f \colon X \to Y$ be a compactifiable morphism of schemes satisfying (\dag) which is flat at dimension $d$. We define
\[ c_{\ell}(X/Y) \coloneqq \mathrm{t}_f \colon \mathbb{Z}_{\ell, X} \to \rdF f^! \mathbb{Z}_{\ell, Y}\dangle{-d} \]
in $\conDb{\etSite{X}, \mathbb{Z}_{\ell, \smdot}}$, i.e., $c_{\ell}(X/Y) \in \mathbb{H}_d(X/Y, \mathbb{Z}_{\ell})$.
\end{Nota}

\begin{Prop}\label{217R0}
Let $X \to S$ and $Y \to S$ be two compactifiable morphisms of schemes satisfying {\upshape(\dag)}, $f \colon X \to Y$ a morphism of $S$-schemes. Assume that $Y \to S$ and $f \colon X \to Y$ are flat at dimension $n$ and $d$ respectively. Then we have
\[ f^{\ast}c_{\ell}(Y/S) = c_{\ell}(X/S) \in \mathbb{H}_{n+d}(X/S, \mathbb{Z}_{\ell}) \,. \]
\end{Prop}

\begin{proof}
This follows from Proposition \ref{195P4} (2).
\end{proof}

\begin{Prop}\label{195P10}
Let $p \colon X \to S$ and $q \colon Y \to S$ be two morphisms of schemes satisfying {\upshape(\dag)} both of which are compactifiable and flat at dimension $d$, $f \colon X \to Y$ a finite $S$-morphism such that $f_{\ast}\mathcal{O}_X$ is a locally free $\mathcal{O}_Y$-module of degree $n$. Then we have
\[ f_{\ast}c_{\ell}(X/S) = n \cdot c_{\ell}(Y/S) \in \mathbb{H}_d(Y, \mathbb{Z}_{\ell}) \,. \]
\end{Prop}

\begin{proof}
By the definition of $f_{\ast}$ and Proposition \ref{195P4} (2), the element $f_{\ast}c_{\ell}(X/S)$ is equal to the composite morphism
\[ \mathbb{Z}_{\ell, Y} \xrightarrow{\delta_f} f_{\ast} \mathbb{Z}_{\ell, X} \xrightarrow{f_{\ast}(\mathrm{t}_f)} f_{\ast} \circ \rdF f^! \mathbb{Z}_{\ell, Y} \xrightarrow{f_{\ast} \circ \rdF f^!(\mathrm{t}_q)} f_{\ast} \circ \rdF f^! \circ \rdF q^! \mathbb{Z}_{\ell, Y} \xrightarrow{\theta_f} \rdF q^! \mathbb{Z}_{\ell, Y} \,. \]
By Diagram \eqref{Eq:4c} and Proposition \ref{195p4a}, we have a commutative diagram
\[ \xymatrix@C+2em@R+1ex{\mathbb{Z}_{\ell, Y} \ar[r]^-{\delta_f} \ar@/_1pc/[drr]^{n \cdot \mathrm{id}} & f_{\ast} \mathbb{Z}_{\ell, X} \ar[r]^-{f_{\ast}(\mathrm{t}_f)} \ar[dr]^{\Tr_f} & f_{\ast} \circ \rdF f^! \mathbb{Z}_{\ell, Y} \ar[r]^-{f_{\ast} \circ \rdF f^!(\mathrm{t}_q)} \ar[d]^{\theta_f} & f_{\ast} \circ \rdF f^! \circ \rdF q^! \mathbb{Z}_{\ell, Y} \ar[d]^{\theta_f} \\
  & & \mathbb{Z}_{\ell, Y}  \ar[r]^-{\mathrm{t}_q} & \rdF q^! \mathbb{Z}_{\ell, Y}} \]
Thus we get the proof.
\end{proof}

\begin{Nota}
Let $X \to S$ be a compactifiable morphism of schemes satisfying (\dag), $i \colon Y \hookrightarrow X$ a closed immersion. Assume that the morphism $Y \to S$ is flat at dimension $d$. Then we define
\[ \widetilde{\clM}_{X/S, \ell}(Y) \coloneqq i_{\ast}c_{\ell}(Y/S) \in \mathbb{H}_d(X/S, \mathbb{Z}_{\ell}) \,. \]
\end{Nota}

\begin{Nota}
Let $X$ be an algebraic scheme over $k$. Then for each $n \in \mathbb{Z}$, there is a canonical homomorphism of groups
\[ \widetilde{\clM}_{X, \ell} \colon \cycZ_n(X) \to \mathbb{H}_n(X, \mathbb{Z}_{\ell}) \,, \qquad \sum a_i \cdot [Y_i] \mapsto \sum a_i \cdot \widetilde{\clM}_{X/k,\ell}(Y_i) \,. \]
\end{Nota}

\begin{Prop}\label{217P1}
Let $f \colon X \to Y$ be a proper morphism of algebraic schemes over $k$. Then for every $n \in \mathbb{N}$, we have a commutative diagram
\[ \xymatrix{\cycZ_n(X) \ar[r]^-{\widetilde{\clM}_{X,\ell}} \ar[d]_{f_{\ast}} & \mathbb{H}_n(X, \mathbb{Z}_{\ell}) \ar[d]^{f_{\ast}} \\
  \cycZ_n(Y) \ar[r]^-{\widetilde{\clM}_{Y,\ell}} & \mathbb{H}_n(Y, \mathbb{Z}_{\ell})} \]
\end{Prop}

\begin{proof}
Let $X'$ be a $n$-dimensional subvariety of $X$, $Y' \coloneqq f(X')$, $i \colon X' \hookrightarrow X$ and $j \colon Y' \hookrightarrow Y$ the inclusion, $g \colon X' \to Y'$ the induced morphism. By Proposition \ref{195P8}, we have
\[ f_{\ast} \circ \widetilde{\clM}_{X,\ell}\bigl([X']\bigr) = f_{\ast} \circ i_{\ast} c_{\ell}(X'/k) = j_{\ast} \circ g_{\ast} c_{\ell}(X'/k) \in \mathbb{H}_n(X, \mathbb{Z}_{\ell}) \,. \]
Since $f_{\ast}[X'] = \deg(X'/Y')[Y']$ (see \cite[1.4]{MR732620}), we have only to prove that
\[ g_{\ast} c_{\ell}(X'/k) = \deg(X'/Y') \cdot c_{\ell}(Y'/k) \in \mathbb{H}_n(Y', \mathbb{Z}_{\ell}) \,. \]

Case 1. $\dim Y' < n$. Then $\deg(X'/Y') = 0$. And by Proposition \ref{195P9}, $\mathbb{H}_n(Y', \mathbb{Z}_{\ell}) = 0$.

Case 2. $\dim Y' = n$. We apply the result in \cite[Ex.~3.7]{MR0463157}. Since the morphism $g$ is generically finite and $Y'$ is an integral scheme, there exists an nonempty subscheme $V$ of $Y'$ such that $g \colon g^{-1}(V) \to V$ is a finite morphism and $g_{\ast}\mathcal{O}_{X'}|_V$ is a locally free $\mathcal{O}_V$-module. Now the proposition follows from Proposition \ref{78fu0i} and Proposition \ref{195P10}.
\end{proof}

\begin{Prop}\label{217P0}
Let $X$ be an algebraic scheme over $k$, $Y$ a $n$-equidimensional closed subscheme of $X$. Then we have
\[ \widetilde{\clM}_{X/k, \ell}(Y) = \widetilde{\clM}_{X,\ell}\bigl([Y]\bigr) \in \mathbb{H}_n(X, \mathbb{Z}_{\ell}) \,. \]
\end{Prop}

\begin{proof}
This is easily deduced from Proposition \ref{3hj01r}.
\end{proof}

\begin{Prop}\label{217U0}
Let $f \colon X \to Y$ be a flat morphism of relative dimension $d$ of algebraic schemes over $k$. Then for every $n \in \mathbb{N}$, we have a commutative diagram
\[ \xymatrix{\cycZ_n(Y) \ar[r]^-{\widetilde{\clM}_{Y,\ell}} \ar[d]^{f^{\ast}} & \mathbb{H}_n(Y, \mathbb{Z}_{\ell}) \ar[d]^{f^{\ast}} \\
  \cycZ_{n+d}(X) \ar[r]^-{\widetilde{\clM}_{X,\ell}} & \mathbb{H}_{n+d}(X, \mathbb{Z}_{\ell})} \]
\end{Prop}

\begin{proof}
Let $\alpha \in \cycZ_n(Y)$. We may assume that $Y$ is a variety of dimension $n$ and $\alpha = [Y]$. Then we have only to apply Proposition \ref{217R0}.
\end{proof}

Now we could prove that $\widetilde{\clM}$ annihilates the rational equivalence.

\begin{Lem}\label{edc1o}
Let $X$ be a nonsingular variety of dimension $n$ over $k$, $D$ an effective divisor on $X$. Then
\[ \widetilde{\clM}_{X/S, \ell}(D) = c_1\bigl(\mathcal{O}(D)\bigr) \in \mathbb{H}_{n-1}(X, \mathbb{Z}_{\ell}) = \mathbb{H}^1(X, \mathbb{Z}_{\ell}) \,. \]
\end{Lem}

\begin{proof}
See \cite[(3.26)]{MR929536}.
\end{proof}

\begin{Thm}
Let $X$ be an algebraic scheme over $k$. Then for each $n \in \mathbb{N}$,
\[ \cycR_n(X) \subseteq \mathrm{Ker}\bigl(\widetilde{\clM}_{X,\ell} \colon \cycZ_n(X) \to \mathbb{H}_n(X, \mathbb{Z}_{\ell})\bigr) \,, \]
i.e., the homomorphism $\widetilde{\clM}_{X,\ell}$ factors through $\cycA_n(X)$. We use $\clM_{X,\ell}$ or $\clM_{\ell}$ or $\clM_X$ to denote the induced homomorphism $\cycA_n(X) \to \mathbb{H}_n(X, \mathbb{Z}_{\ell})$.
\end{Thm}

\begin{proof}
After applying \cite[Proposition 1.6]{MR732620} together with Proposition \ref{217P1} and Proposition \ref{217U0}, we have only to prove that
\[ \widetilde{\clM}_{\mathbf{P}^1_k, \ell}(0) = \widetilde{\clM}_{\mathbf{P}^1_k, \ell}(\infty) \in \mathbb{H}_0(\mathbf{P}^1_k, \mathbb{Z}_{\ell}) \,. \]
This is by Lemma \ref{edc1o}.
\end{proof}

In the following, we define the degree map for the homology of degree zero. Note that $\mathbb{H}_0(\Spec k, \mathbb{Z}_{\ell}) = \mathbb{Z}_{\ell}$. So we have

\begin{Defi}
 For any proper algebraic scheme $X$ over $k$, we define \emph{degree map} $\deg_{\ell}$ to be the homomorphism
\[ \mathbb{H}_0(X, \mathbb{Z}_{\ell}) \xrightarrow{p_{\ast}} \mathbb{H}_0(\Spec k, \mathbb{Z}_{\ell}) = \mathbb{Z}_{\ell} \,, \]
where $p \colon X \to \Spec k$ is the structural morphism.
\end{Defi}

\begin{Lem}\label{217L0}
Let $X$ be a $n$-dimensional proper algebraic scheme over $k$.
\begin{enumerate}
\item We have a commutative diagram.
\[ \xymatrix{\upH{2n}\bigl(\etSite{X}, \mathbb{Z}_{\ell}(n)\bigr) \ar[rr]^-{(\mathrm{t}_X)_{\ast}} \ar[rd]_-{\Tr_X} & & \mathbb{H}_0(X, \mathbb{Z}_{\ell}) \ar[ld]^{\deg_{\ell}} \\ & \mathbb{Z}_{\ell}} \]
\item We have a commutative diagram.
\[ \xymatrix@C+0.5em{\cycA_0(X) \ar[r]^-{\clM_{X, \ell}} \ar[d]_{\deg} & \mathbb{H}_0(X, \mathbb{Z}_{\ell}) \ar[d]^{\deg_{\ell}} \\ \mathbb{Z} \ar@{^(->}[r] & \mathbb{Z}_{\ell}} \]
\end{enumerate}
\end{Lem}

\begin{proof}
(1) is by the commutative diagram \ref{Eq:4c}.

(2) is by Proposition \ref{217P1}.
\end{proof}

\begin{Prop}
Assume that $k$ is separably closed and let $X$ be a nonsingular complete variety over $k$. Then $\deg_{\ell} \colon \mathbb{H}_0(X, \mathbb{Z}_{\ell}) \isoTo \mathbb{Z}_{\ell}$ is an isomorphism.
\end{Prop}

\begin{proof}
Put $\dim X = n$. By Proposition \ref{195T5}, we have only to prove that $\Tr_X \colon \upH{2n}\bigl(\etSite{X}, \mathbb{Z}_{\ell}(n)\bigr) \to \mathbb{Z}_{\ell}$ is an isomorphism. This is by \cite[VI, 11.1 (a)]{MR559531}.
\end{proof}

The following theorem shows that the cycle map $\clM_{X, \ell}$ annihilate algebraic equivalence of cycles.

\begin{Thm}
Assume that $k$ is separably closed and let $X$ be an algebraic scheme over $k$. Then for each $n \in \mathbb{N}$, the cycles in $\cycA_n(X)$ which are algebraically equivalent to zero (in the sense of \cite[10.3]{MR732620}), are contained in $\mathrm{Ker}\bigl(\clM_{X, \ell})$.
\end{Thm}

\begin{proof}
By Proposition \ref{base0} we may assume that $k$ is algebraically closed. Let $c_1, c_2 \in \cycA_n(X)$ such that $c_1 \sim_{\textup{a}} c_2$; and let $T$ be a nonsingular curve over $k$. $t_1, t_2 \in T(k)$, $c \in \cycA_{n+1}(X \times_k T)$ such that $c_{t_i} = c_i$ for $i=1,2$. Obviously we may assume that $c = [Y]$, where $Y$ is a $(n+1)$-dimensional subvariety of $X \times_k T$ such that for all $t \in T(k)$, $Y$ is not contained in
\[ \iDe \times t_i \colon X \hookrightarrow X \times_k T \,. \]
Obviously the induced morphism $Y \to T$ is dominant and flat. Put
\[ \alpha \coloneqq \widetilde{\clM}_{(X \times_k T)/T, \, \ell}(Y) \in \mathbb{H}_n\bigl((X \times_k T)/T, \mathbb{Z}_{\ell}\bigr) \,. \]
By Proposition \ref{j4i0yn} and Proposition \ref{195P11} (2), we have
\[ t_i^{\ast}(\alpha) = \widetilde{\clM}_{X, \ell}(Y_{t_i}) = \clM_{X, \ell}(c_{t_i}) \,. \]
So we have only to apply Proposition \ref{eaec8}.
\end{proof}

\section{Cap-products and Compatibility with Chern classes}\label{Sec8}

First we define the cap-products for the $\ell$-adic homology.

\begin{Defi}[Cap-Product]\label{j75ol}
Let $f \colon X \to Y$ and $g \colon Y \to Z$ be compactifiable morphisms of schemes satisfying (\dag), $\mathscr{M}$ and $\mathscr{N}$ two objects in $\conDn{\etSite{Z}, R_{\smdot}}$.

For every $m, n \in \mathbb{Z}$, there is a \emph{cap-product}
\[ \rmH{m}(X \xrightarrow{f} Y, g^{\ast}\mathscr{M}) \times \rmH{n}(Y \xrightarrow{g} Z, \mathscr{N}) \xrightarrow{\cap} \rmH{m+n}(X \xrightarrow{g \circ f} Z, \mathscr{M} \ttP[R] \mathscr{N}) \,, \]
defined as follows. Let $\alpha \in \rmH{m}(X \xrightarrow{f} Y, g^{\ast}\mathscr{M})$ and $\beta \in \rmH{n}(Y \xrightarrow{g} Z, \mathscr{N})$, then we define $\alpha \cap \beta$ to be the composite morphism
\[ R_X \xrightarrow{\alpha \ttP f^{\ast}\beta} \rdF f^! \circ g^{\ast} \mathscr{M}[-m] \ttP[R] f^{\ast} \circ \rdF g^! \mathscr{N}[-n] \xrightarrow{\varphi} \rdF (g \circ f)^!(\mathscr{M} \ttP[R] \mathscr{N})[-(m+n)] \,. \]
where $\varphi$ is defined in Proposition \ref{195P13}.

Similarly we may define the \emph{cap-product} for $\mathbb{H}_{\ast}$ as follows:
\[ \mathbb{H}_m(X \xrightarrow{f} Y, g^{\ast}\mathscr{M}) \times \mathbb{H}_n(Y \xrightarrow{g} Z, \mathscr{N}) \xrightarrow{\cap} \mathbb{H}_{m+n}(X \xrightarrow{g \circ f} Z, \mathscr{M} \ttP[R] \mathscr{N}) \,, \]
\end{Defi}

In particular if $X \to S$ is a compactifiable morphisms of schemes satisfying (\dag), and $\mathscr{N}$ an object in $\conDn{\etSite{S}, R_{\smdot}}$, then for every $m, n \in \mathbb{Z}$, there are cap-products
\begin{gather*}
	\upH{m}(X, R) \times \rmH{n}(X/S, \mathscr{N}) \xrightarrow{\cap} \rmH{n-m}(X/S, \mathscr{N}) \,, \\
	\mathbb{H}^m(X, R) \times \mathbb{H}_n(X/S, \mathscr{N}) \xrightarrow{\cap} \mathbb{H}_{n-m}(X/S, \mathscr{N}) \,.
\end{gather*}

The following Proposition can be directly calculated.

\begin{Prop}[Projection Formula]\label{p7u8i0}
	Let $f \colon X \to Y$ and $g \colon Y \to S$ be morphisms of schemes satisfying {\upshape(\dag)} with $f$ proper and $g$ compactifiable, $\mathscr{N}$ an object in $\conD{\etSite{Y}, R_{\smdot}}$. Then we have
	\begin{enumerate}
		\item For every $\alpha \in \upH{r}(Y, R)$ and $\beta \in \rmH{n}(X/S, \mathscr{N})$, we have
			\[ \alpha \cap f_{\ast}(\beta) = f_{\ast}\bigl(f^{\ast}(\alpha) \cap \beta\bigr) \in \rmH{n-r}(Y/S, \mathscr{N}) \,. \]
		\item For every $\alpha \in \mathbb{H}^r(Y, R)$ and $\beta \in \mathbb{H}_n(X/S, \mathscr{N})$, we have
			\[ \alpha \cap f_{\ast}(\beta) = f_{\ast}\bigl(f^{\ast}(\alpha) \cap \beta\bigr) \in \mathbb{H}_{n-r}(Y/S, \mathscr{N}) \,. \]
	\end{enumerate}
\end{Prop}

It may be further showed that the cup-product defined in Definition \ref{j75ol} has many similar properties with bivariant intersection theory defined in \cite[Ch.~17]{MR732620}, i.e., has associativity and is compatible with the Pull-back functor $f_{\ast}$, the push-out functor $f^{\ast}$ and the base change functor $u^{\ast}$. Since we need not them here, so we leave it to the readers.

\medskip

Next, we review the cycle maps for locally free sheaves. First by \cite[(3.26) a)]{MR929536}, we have a homomorphism of groups
\begin{equation}\label{E:clinv}
c_1^{\ell} \colon \Pic X \to \mathbb{H}^1(X, \mathbb{Z}_{\ell})
\end{equation}
for every scheme $X$ satisfying (\dag). The following two propositions depict the cycle maps for locally free sheaves.

\begin{Prop}
Let $S$ be a scheme satisfying {\upshape(\dag)}, $\mathscr{E}$ a locally free $\mathcal{O}_S$-module of constant rank $r+1$, $P \coloneqq \aProj(\mathscr{E})$, $p \colon P \to S$ the projection. Then  there is a canonical isomorphism of $\mathbb{Z}_{\ell}$-algebras
\[ \mathbb{H}^{\ast}(S, \mathbb{Z}_{\ell})[T]/(T^{r+1}) \isoTo \mathbb{H}^{\ast}(P, \mathbb{Z}_{\ell}) \,, \qquad \overline{T} \mapsto c_1^{\ell}\bigl(\mathcal{O}_P(1)\bigr) \,. \]
\end{Prop}

\begin{proof}
See \cite[(6.13)]{MR929536} or \cite[VII, 2.2.6]{SGA5}.
\end{proof}

\begin{Prop}\label{195R3}
Let $X$ be a scheme satisfying {\upshape(\dag)}, $\mathscr{E}$ a locally free $\mathcal{O}_S$-module of constant rank $m$, $P \coloneqq \aProj(\mathscr{E}\spcheck)$, $p \colon P \to X$ the projection, $\xi \coloneqq c_1^{\ell}\bigl(\mathcal{O}_P(1)\bigr)$. Then there exists a unique element $c_r^{\ell}(\mathscr{E}) \in \mathbb{H}^r(X, \mathbb{Z}_{\ell})$ for each $r \in \mathbb{N}$, such that
\[ \left\{\begin{array}{l}
  \sum\limits^m_{i=0} c_i^{\ell}(\mathscr{E}) \xi^{m-i} = 0, \\[2ex]
  c_0^{\ell}(\mathscr{E}) = 1, \\[.5ex]
  c_r^{\ell}(\mathscr{E}) = 0 \text{ for } r > m.
\end{array}\right. \]
\end{Prop}

Now we define the trace morphisms for regular immersions of codimension $1$. Let $X$ be a scheme satisfying (\dag) and $i \colon D \hookrightarrow X$ a regular closed immersion of codimension $1$. By \cite[(3.26) and the proof]{MR929536}, $i \colon D \hookrightarrow X$ determinates an element
\[ \mathrm{t}_i \in \mathrm{H}^2_{D, \mathrm{cont}}\bigl(\etSite{X}, \mathbb{Z}_{\ell}(1)\bigr) = \Hom_{\conD{\etSite{D}, \mathbb{Z}_{\ell, \smdot}}}\bigl(\mathbb{Z}_{\ell}, \rdF i^! \mathbb{Z}_{\ell}\dangle{1}\bigr) \,. \]
Similar to \cite[(cycle) 2.3.1]{MR0463174}, we have

\begin{Prop}\label{rei08u1}
Let $S$ be a scheme satisfying {\upshape(\dag)}, $f \colon X \to S$ and $g \colon Y \to S$ two compactifiable morphisms which are flat at dimension $n$ and $n-1$ respectively, $i \colon Y \hookrightarrow X$ a regular closed immersion of codimension $1$ such that $f \circ i = g$. Then we have
\begin{enumerate}
\item The composite morphism
\[ \rdF g_! \mathbb{Z}_{\ell} \xrightarrow{\rdF g_!(\mathrm{t}_i)} \rdF g_! \circ \rdF i^!\mathbb{Z}_{\ell}\dangle{1} = \rdF f_! \circ \rdF i_{\ast} \circ \rdF i^! \mathbb{Z}_{\ell}\dangle{1} \xrightarrow{\rdF f_!(\theta_i)} \rdF f_! \mathbb{Z}_{\ell}\dangle{1} \xrightarrow{\Tr_f} \mathbb{Z}_{\ell}\dangle{-(n-1)} \]
is equal to $\Tr_g$.
\item The composite morphism
\[ \mathbb{Z}_{\ell} \xrightarrow{\mathrm{t}_i} \rdF i^! \mathbb{Z}_{\ell}\dangle{1} \xrightarrow{\rdF i^!(\mathrm{t}_f)} \rdF i^! \circ \rdF f^! \mathbb{Z}_{\ell}\dangle{-(n-1)} = \rdF g^! \mathbb{Z}_{\ell}\dangle{-(n-1)} \]
is equal to $\mathrm{t}_g$.
\end{enumerate}
\end{Prop}

Finally we could prove that the cycle maps are compatible with Chern classes. According to \cite[Ch.~3]{MR732620}, if $X$ is an algebraic scheme over $k$ and $\mathscr{E}$ is a locally free $\mathcal{O}_X$-module, then there is an operation of Chern classes on each Chow group
\[ \cycA_r(X) \to \cycA_{r-i}(X) \,, \qquad \alpha \mapsto c_i(\mathscr{E}) \cap \alpha \,. \]

\begin{Thm}
Let $X$ be an algebraic scheme over $k$, $\mathscr{E}$ a locally free $\mathcal{O}_X$-module, $\alpha \in \cycA_r(X)$. Then we have
\begin{equation}\label{E:fjk78u}
c_i^{\ell}(\mathscr{E}) \cap \clM_{X,\ell}(\alpha) = \clM_{X,\ell}\bigl(c_i(\mathscr{E}) \cap \alpha) \in \mathbb{H}_{r-i}(X, \mathbb{Z}_{\ell}) \,.
\end{equation}
\end{Thm}

\begin{proof}
By the the projection formulas (Proposition \ref{p7u8i0} and \cite[Theorem 3.2 (c)]{MR732620}), we obtain that if $f \colon X' \to X$ is a proper morphism and $\alpha' \in \cycA_r(X')$ such that $f_{\ast}(\alpha') = \alpha$ and the pair $(f^{\ast}\mathscr{E}, \alpha')$ satisfies \eqref{E:fjk78u}, then the pair $(\mathscr{E}, \alpha)$ also satisfies \eqref{E:fjk78u}. Thus by the splitting construction (see \cite[\S3.2]{MR732620}), we may assume that $\mathscr{E} = \mathscr{L}$ is an invertible $\mathcal{O}_X$-module and have only to prove that
\begin{equation}\label{E:fjk78v}
c_1^{\ell}(\mathscr{L}) \cap \clM_{X,\ell}(\alpha) = \clM_{X,\ell}\bigl(c_1(\mathscr{L}) \cap \alpha\bigr) \in \mathbb{H}_{r-1}(X, \mathbb{Z}_{\ell}) \,.
\end{equation}
Moreover we may assume that $X$ is a variety of dimension $r$ and $\alpha = [X]$. After replacing $X$ with its normalization, we may assume that $X$ is normal. Then we may assume that $\mathscr{L} = \mathcal{O}(Y)$ where $Y \hookrightarrow X$ is a regular closed immersion of codimension $1$. Then we have only to apply Proposition \ref{rei08u1}.
\end{proof}

\bibliographystyle{plain}
\bibliography{Cycle}

\end{document}